\documentclass[12pt]{article}

\usepackage[T1]{fontenc}
\usepackage[utf8]{inputenc}
\usepackage[frenchb]{babel}
\usepackage[bookmarks,hidelinks,pdftex]{hyperref}
\usepackage[usenames, dvipsnames]{xcolor}
\usepackage{graphicx}
\DeclareGraphicsExtensions{.jpg}

\usepackage{tikz}
\title{Paul Mansion (1844-1919) : more than 400 academic  publications\\
Centenary Paul Mansion ({\it working paper, 2019})\\
Paul Mansion (1844-1919) : plus de 400 publications mathématiques selon le Jahrbuch\\
Centenaire Paul Mansion ({\it document de travail, 2019})}
\author{Hervé Le Ferrand\footnote{Institut de Mathématiques de Bourgogne, leferran@u-bourgogne.fr}}
\date{\today}

\begin{document}

\maketitle

According to  the Liber Memorialis of the University of Ghent, the Belgian mathematician Paul Mansion (1844-1919) has published more than 349 academic papers and books. For our part, we were able to calculate the correct number by using the journal {\it Das Jahrbuch \"uber die Fortschritte der Mathematik} (1869-1942). We concluded that Paul Mansion has published about  400 academic papers and books.

\section{Introduction}
De sa nomination en 1867 à l'Université de Gand jusqu'à la fin de sa vie, Paul Mansion\footnote{Pour des éléments sur la vie et l'oeuvre de Paul Mansion, on peut consulter \cite{lef3}.} a publié un nombre impressionnant d'articles scientifiques. On retrouve cette même force de travail dans ses activités éditoriales et dans sa correspondance scientifique. Il suffit par exemple de consulter le fonds Paul Mansion conservé à la bibliothèque royale de Belgique à Bruxelles, pour se rendre compte du rayonnement scientifique du mathématicien de l'Université de Gand. 

La biographie que l'on peut lire sur le site Mac Tutor\footnote{Mac Tutor History of Mathematics archive : \\http://www-history.mcs.st-and.ac.uk/Biographies/Mansion.html} mentionne :
\begin{quote}
\begin{it}
Mansion was highly productive. When Alphonse Demoulin wrote the obituary of Mansion [3] (published in 1929) he included a list of 349 of his works although in fact the list contains many papers with the same title published in parts over many years which each appear as a single item. Nor were these works merely of local interest, for many were considered important enough to be translated into German or to be republished in foreign publications (...)
\end{it}
\end{quote}
Combien d'articles scientifiques Paul Mansion a-t-il  publié ? Doit-on conserver ce nombre de 349 ? Nous allons, pour répondre à cela, comparer la liste établie par Paul Mansion lui-même avec celle que l'on peut composer à partir du {\it Jahrbuch}.

\section{Le nombre de publications mathématiques de Paul Mansion selon le Jahrbuch}
En 1927 devant la Classe des Sciences de l'Académie royale des sciences, des lettres et des beaux-arts de Belgique\footnote{Cette académie est divisée en quatre classes.}, Alphonse Demoulin\footnote{Alphonse Demoulin est  un ancien élève et collègue de Paul Mansion. Alphonse Demoulin étudia et enseigna à l'Université de Gand à partir de 1893 (jusqu'en 1936)\cite{doc2}. Il fut aussi l'élève de Gaston Darboux (1842-1917) à Paris en 1892. Demoulin travailla essentiellement dans le domaine de la Géométrie Différentielle et ses travaux furent récompensés par troix prix de l'Académie des Sciences de Paris : le prix du baron de Joest en 1906, le prix Bordin en 1911\footnote{Sujet proposé : \og {\it perfectionner en un point important la théorie des systèmes triples de surfaces orthogonales} \fg.} et le prix Poncelet en 1945.}\footnote{On pourra aussi consulter \cite{lef2}.}(1869-1947), prononce un éloge de Paul Mansion \cite{Demoulin}. Dans le texte de ce discours \cite{Demoulin}, Alphonse Demoulin donne la liste des publications de Paul Mansion telle que ce dernier l'a établie en 1913 dans le Liber Memorialis de l'Université de Gand \cite{doc1}. On y trouve à la fois une liste de 349 entrées rangées par titres de journaux et d'éditeurs pour les livres,  et une classification thématique des écrits de Paul  Mansion selon les domaines mathématiques. D'un autre côté, la version électronique du {\it Das Jahrbuch \"uber die Fortschritte der Mathematik}, The Jahrbuch Project Electronic Research Archive for Mathematics (ERAM)\footnote{https://www.emis.de/MATH/JFM/} donne 478 entrées\footnote{Avec quelques doublons.} lorsqu'on entre le nom \og Paul Mansion\fg. Rappelons que le {\it Jahrbuch} recense les articles de mathématiques parus de 1868 à 1942. Le document \cite{doc3} de l'Université de Heidelberg donne des références historiques sur ce journal allemand. On apprend par exemple que Paul Mansion a été référé du {\it Jahrbuch} de 1871 à 1916.

Dans sa propre liste \cite{doc1} parue dans le Liber Memorialis, Paul Mansion a effectué des regroupements. Ainsi un item peut correspondre à plusieurs articles. Pour illustration, l'entrée portant le numéro 72 (revue {\it Mathesis}) intitulé {\it Théorie des limites et des infiniment petits} correspond à 8 articles parus dans différents numéros de {\it Mathesis}. 

Dans l'ordre, Paul Mansion a indiqué les revues suivantes\footnote{Nous indiquons, si besoin, entre parenthèses les initiales utilisées par le {\it Jahrbuch}.} :
\begin{itemize}
\begin{it}
\item Bulletins de Belgique (Bull. de Belg.)\
\item Bulletins de la Classes des Sciences (Belge Bull. Sciences)
\item Annuaire (Annuaire Ac. de Belgique)
\item Biographie nationale
\item Revue de l'Instruction publique en Belgique
\item Nouvelle correspondance mathématique (N. C. M.)
\item Mathesis (Math.){Revue fondée par Paul Mansion et Neuberg.}
\item Messenger of Mathematics (Messenger ou Mess.)
\item Report of the British Association for the Advance of Sciences (Rep. Brit. Ass.)
\item Archives de Grunert, Archiv der Mathematik und Physik
\item Bull. di Bibliografia e di Storia delle Scienze matem. e fisiche (Boncompagni\footnote{Balthazar Boncompagni (1821-1894)} Bull.)
\item Comptes rendus des séances de l'Académie des Sciences de Paris (C.R.)
\item Annales de la Société scientifique de Bruxelles (Ann. Soc. scient. Brux.)
\item Revue des questions scientifiques (Rev. d. qu. sc.)
\item Mémoire de la Société royale des Sciences de Liège
\item Bibliotheca mathematica de M. G. Enestr\"om\footnote{Gustaf Enestr\"om (1852-1923)} (Biblioth. Mathem.)
\item Bulletin des sciences mathématiques ( Bulletin de Darboux\footnote{Gaston Darboux (1842-1917)}, Darb. Bull.)
\item Wiadomosci Matematyszne, de Dickstein (Wiad. mat.)
\item Congrès scientifique international des catholiques
\item Revue néo-scolastique
\item Jahresbericht der Deutschen Mathematiker-Vereinigung (Deutsche Math. Ver.)
\item Abhandlungen zu Geschichte der Mathematik
\item Enseignement Mathématique (Ens. Math.)
\item Le Mathematiche pure ed applicate (Mat. pure ed appl.)
\item Congrès de l'enseignement moyen
\item Revue des humanités en Belgique
\item Acta Mathematica (Acta Math.)
\item Revue de philosophie de Peilhaube
\item Bulletin de la Société physico-mathématique de Kazan
\item Congrès international d'expansion économique mondiale
\item Moniteur belge
\end{it}
\end{itemize}
auxquelles s'ajoutent des \og écrits publiés dans divers recueils et dans divers ouvrages\fg, des \og ouvrages publiés à part\fg\ et des \og traductions\fg. 

Donnons quelques détail. Les {\it Bulletins de Belgique} sont les {\it Bulletins de l'Académie royale des sciences, des lettres et des beaux-arts de Belgique} publiés de 1846 à 1898. En 1899 apparaît le {\it 
Bulletin de la Classe des sciences (Académie royale de Belgique)}. La {\it Nouvelle correspondance mathématique} est fondée par Eugène Catalan (1814-1894) et Paul Mansion en 1874. La revue {\it Mathesis} est créée par Paul Mansion et Joseph Neuberg (1840-1926) en 1881 et remplace la  {\it Nouvelle correspondance mathématique}. Le {\it Messenger of Mathematics}, journal anglo-irlandais, paraît de 1872 à 1929. La revue {\it Archiv der Mathematik und Physik} est publiée de 1841 à 1920. Boncompagni édite  le {\it Bullettino di bibliografia e di storia delle scienze matematiche e fisiche} de 1868 à 1887. La {\it Revue des questions scientifiques} est éditée par la Société Scientifique de Bruxelles à partir de 1877\footnote{Après une interruption de plusieurs années, elle est de nouveau publiée depuis deux ans : https://www.rqs.be/app/views/index.php}. De 1884 à 1914 Gustaf Enestr\"om est l'éditeur de {\it Bibliotheca mathematica}. La revue polonaise {\it Wiadomosci Matematyszne} paraît de 1897 à 1939. Le journal allemand {\it Jahresbericht der Deutschen Mathematiker-Vereinigung} est fondé en 1890 par les mathématiciens allemands Georg Cantor (1845-1918), Walther van Dyck (1856-1934) et Emil Lampe (1840-1918). En 1877 paraît le premier volume du {\it Abhandlungen zur Geschichte der mathematischen Wissenschaften}. Le mathématicien suisse Henri Fehr (1870-1954) et le mathématicien et homme politique français Charles-Ange Laisant(1841-1920) fondent  l'{\it Enseignement mathématique} en 1899 . {\it Le Mathematiche pure ed applicate} est fondé par Cristoforo Alasia De Quesada (1869 -1918) en 1901\footnote{Voir : https://cirmath.hypotheses.org/files/2015/09/Colloque-Cirmath-IML-Maria-Rosaria-Enea.pdf}. 

On peut remarquer que Paul Mansion a écrit en anglais, allemand, néerlandais et vraisemblablement en italien. Pour certains journaux de la liste ci-dessus, les articles publiés dans ces journaux ne sont pas référencés par le Jahrbuch. C'est le cas notamment de {\it Congrès scientifique international des catholiques}, {\it Revue néo-scolastique} ou encore {\it Revue de philosophie de Peilhaube} même si les articles de Paul Mansion portent sur les mathématiques.

Revenons au {\it Jahrbuch}. La liste que nous avons pu établir à partir de la version électonique du {\it Jahrbuch} contient, comme nous l'avons déjà indiqué, 478 items. Il y en a réellement 472 différents car quelques publications ont été répertoriées plusieurs fois. On peut considérer que le nombre de publications mathématiques de Paul Mansion dépasse 400, ordre de grandeur de la liste du {\it Jahrbuch}. Ceci est cohérent avec la liste établie par Paul Mansion dans \cite{doc1} car il faut tenir compte des regroupements effectués par Paul Mansion.

Qu'en est-il du \og rythme\fg\  de publication de Paul Mansion ? Si on classe par années les items de la liste du {\it Jahrbuch}, on obtient le tableau suivant :

\begin{table}[h]

\begin{center}
\begin{tabular}{|c|c|}
\hline
années&nombres de publications\\
\hline
1868-1880&103\\
\hline
1881-1890&106\\
\hline
1891-1900&125\\
\hline
1901-1910&89\\
\hline
1911-1922&49\\
\hline
\end{tabular}
\end{center}
\caption{publications par années}
\end{table}
On note une activité de publication très soutenue jusqu'en 1910 date de son éméritat.

\section{Conclusion}
Il nous a paru important, cent ans après la disparition de Paul Mansion, de faire un éclairage sur la liste de ses publications  pour en montrer d'une part la richesse et d'autre part pour rendre hommage  à ce mathématicien dont le rayonnement scientifique ne peut qu'impressionner. 

Ce présent document est un document de travail qui est appelé à être complété.

\end{document}